\documentclass[reqno]{amsart}

\usepackage{xcolor}
\usepackage[nobysame,initials]{amsrefs}
\usepackage{hyperref}
\usepackage{orcidlink}
\usepackage{float}

\usepackage{bbm}
\usepackage{todonotes}
\usepackage{cancel}

\usepackage{pgf}
\usepackage{tikz}
\usepackage{pgfplots}
\usetikzlibrary{arrows}
\usetikzlibrary{intersections,calc}

\allowdisplaybreaks

\newtheorem{theorem}{Theorem}[section]
\newtheorem{corollary}[theorem]{Corollary}
\newtheorem{lemma}[theorem]{Lemma}

\numberwithin{equation}{section}

\newcommand{\R}{\mathbb{R}}
\newcommand{\EE}{\mathbb{E}}
\newcommand{\dx}{{\rm d}}

\DeclareMathOperator{\cl}{cl}
\DeclareMathOperator{\vol}{Vol}

\title[On generalized disc-polygons in plane convex bodies]{On generalized disc-polygons in plane convex bodies with a higher degree of smoothness}

\author[F.~Fodor]{Ferenc Fodor$^{\orcidlink{0000-0001-9747-1981}}$} 
\address{Bolyai Institute, University of Szeged, 
Aradi v\'ertan\'uk tere 1, H-6720 Szeged, Hungary} 
\email{fodorf@math.u-szeged.hu}

\author[D.~I.~Papv\'ari]{D\'aniel I. Papv\'ari$^{\orcidlink{0009-0008-4813-2906}}$}
\address{Bolyai Institute, University of Szeged, 
	Aradi v\'ertan\'uk tere 1, H-6720 Szeged, Hungary}
\email{papvari.daniel.istvan@stud.u-szeged.hu}


\pgfplotsset{compat=1.18}

\begin{document}
    \begin{abstract}
        We prove power series expansions for the expectations of the number of vertices and missed area of random $L$-convex polygons in planar convex bodies with sufficiently smooth boundaries. Random $L$-convex polygons arise as the intersection of all translates of a fixed convex set $L$ that contain i.i.d. uniform random points from a suitable plane convex body $K$. Our results extend the asymptotic formulas proved in Fodor, Papv\'ari and V\'\i gh \cite{FPV20} and Fodor and Montenegro \cite{FM24}, and have consequences about $L$-convex floating bodies and relative affine surface area that were investigated by Sch\"utt, Werner and Yalikun in \cite{SWY25}. 
    \end{abstract}

    \maketitle

\section{Introduction and results}
In random polytope theory, perhaps the most commonly investigated model is when one selects $n$ identically distributed independent (i.i.d.) random points in a convex body (a compact convex set with a non-empty interior) according to a given probability distribution and then takes the convex hull of the points, thereby obtaining a random polytope. The main objects of study are the properties of geometric functionals associated with such random polytopes, for example, intrinsic volumes, the number of $k$-dimensional faces, etc. For the most recent overview of this subject area, see, for example, the survey by Schneider \cite{Sch17} and the numerous references therein; for earlier surveys, see B\'ar\'any \cite{Ba08}, Reitzner \cite{R10}, Weil, and Wieacker \cite{WeWi93}.

Recently, several models have emerged in which the method of generating the random polytope differs from that of taking the classical convex hull. One such prominent example uses the concept of the so-called $L$-convex hull, in which one takes the intersection of all translates of a suitable fixed convex body $L$ containing the random points. Results in this model include asymptotic formulas for expectations, bounds on variances, and even central limit theorems under various geometric conditions. For a comprehensive survey of the state-of-the-art on this topic, see the paper \cite{BLN25+} by Bezdek, L\'angi, and Nasz\'odi. In this paper, we prove series expansions on the expected area and the number of vertices in the planar version of this model under stronger smoothness assumptions. Our results add to those of Fodor and Montenegro \cite{FM24}. 

Our setting is the $d$-dimensional Euclidean space $\R^d$. We denote the Euclidean inner product in $\R^d$ by $\langle\cdot, \cdot\rangle$ and the corresponding norm by $\|\cdot\|$. The origin-centered unit radius ball is $B^d$, and its boundary $S^{d-1}=\partial B^d$ is the unit sphere.
Volume in $\R^d$ (Lebesgue measure) is denoted by $\vol_d(\cdot)$, and in the special case when $d=2$ by $A(\cdot)$. Euler's gamma function is denoted by $\Gamma(\cdot)$.
A convex body in $\R^d$ is a  compact convex set $K$ with a non-empty interior. For information about the properties of convex bodies, we refer to Schneider's monograph \cite{Sch14}.

Next, we introduce our setting by defining the notion of $L$-convex bodies in $\R^d$. This is one of the recently emerged generalizations of classical convexity that appears naturally in many problems. 

Let $L$ be a convex body in $\R^d$, and let $K$ be a compact set that is contained in a translated copy of $L$. We call $K\subset\R^d$ $L$-convex if it is equal to the intersection of all translates of $L$ that contain $K$. A compact set with a non-empty interior that is $L$-convex is called an $L$-convex body (also called an $L$-ball convex body, $L$-strongly convex body, etc.). $L$-convex bodies are also convex in the classical sense. 
In the notion of $L$-convexity, the translates of $L$ play a similar role to closed half-spaces in classical convexity. 

If a set $X\subset \R^d$ is contained in a translate of $L$, then we call the intersection of all translates of $L$ that contain $K$ the $L$-convex hull of $K$, which we denote by $[X]_L$.
If $X$ is not contained in any translate of $L$, then we define $[X]_L=\R^d$; this case will not appear in this paper. 
If $X$ is finite, then $[X]_L$ is called an $L$-polytope and, in $\R^2$, an $L$-polygon.
If $K$ is an $L$-convex body and $X\subset K$, then $[X]_L\subset K$.
For a discussion of the geometric properties of $L$-convex sets, see, for example, L\'angi, Nasz\'odi, Talata \cite{LNT13}, Bezdek et al. \cite{bez05}, and the references in Bezdek, L\'angi and Naszódi \cite{BLN25+}.
  
The $L$-convexity of $K$ gives rise to approximations by random (and non-random) $L$-polytopes in the following way. Let $X_n=\{p_1,\ldots, p_n\}$ be i.i.d. random points from $K$ selected according to the uniform probability distribution (normalized Lebesgue measure). Then $K_n^L=[X_n]_L$ is a (uniform) random $L$-polytope in $K$. In this paper, we study the asymptotic behavior of the expectation of the area difference of $K$ and $K_n^L$ (also called the missed area), that is, $\EE(A(K\setminus K_n^L))$, and the expectation of the number of vertices $\EE (f_0(K_n^L))$. We note that most of the approximation results in this model are $2$-dimensional (including the sphere and the hyperbolic plane) with some notable exceptions, for example, by Marynych and Molchanov \cite{MM22}. This phenomenon is due to the intrinsic geometry of the model, which poses several obstacles that prevent the use of several successful methods from classical random polytope theory.  

Our motivations for proving series expansions of expectations are twofold. One reason to investigate series expansions under stronger differentiability conditions goes back to the results of Gruber \cite{G96} and Reitzner \cites{R01, R04}. Gruber \cite{G96} proved series expansions for the support function and the mean width of uniform random polytopes. Reitzner extended these results to the expectation of all intrinsic volumes and the number of vertices for $d=2$ in \cite{R01} using floating bodies, and for general $d$ in \cite{R04}. Fodor and Montenegro \cite{FM24} proved series expansions for the number of vertices and the missed area of uniform random $L$-polygons in the case $L=RB^2$ for $d=2$. Our results directly extend those of \cite{FM24}.

Our other reason is that in classical approximation results (both random and best), two geometric notions often appear: the convex floating body and the affine surface area. For a convex body $K$, and for sufficiently small $\delta>0$, the classical convex floating body $K_{\delta}$ is the intersection of all closed half-spaces that cut off at most volume $\delta$ from $K$. The difference set $K(\delta)=K\setminus K_\delta$ is referred to as the wet part of $K$ with parameter $\delta$. In random approximation theory, one of the most fundamental results is due to B\'ar\'any and Larman \cite{BL89}, which states that the order of magnitude of the volume of the missed part of a classical random polytope in a convex body $K$ is the same as that of the volume of $K(1/n)$. This statement holds regardless of the smoothness of the boundary of $K$. The connection between the floating body and expectations goes even deeper by, for example, the results of Buchta and Reitzner \cite{BR97}. Reitzner \cite{R01} first proved series expansions for the volume of the wet part $K(\delta)$ under stronger differentiability conditions, then used the result from \cite{BR97} to obtain series expansions for the expectation of the missed area and the number of vertices of random polygons. 

The notion of $L$-convex floating body $K^L_\delta$ and wet part $K^L(\delta)$ may be introduced similarly as the intersection of all translates of $L$ that cut off volume at most $\delta$ from $K$, see, for example, Fodor and Papv\'ari \cite{FP25} for the case when $L=RB^2$, and Sch\"utt, Werner and Yalikun \cite{SWY25} for the general case. 
There is evidence indicating that it is very likely that the order of magnitude of the expectation of the missed volume of uniform random $L$-polytopes is the same as the order of magnitude of the $L$-convex wet part with parameter $1/n$; however, this has not yet been proven. Recent results about the $L$-convex analog of the affine surface area by Sch\"utt, Werner and Yalikun in \cite{SWY25} also support this conjecture.

The affine surface area $as(K)$ is the right derivative of the volume of the convex floating body $K_\delta$ with respect to $\delta$, see, for example, the paper \cite{SW22} by Sch\"utt and Werner. Thus, the affine surface area is naturally present in many (best and random) approximation results, and understanding its properties is fundamental in polytopal approximation and, in general, convex geometry. 

In \cite{SWY25}, Sch\"utt, Werner and Yalikun investigate the $L$-convex floating body and introduce the relative affine surface area $as^L(K)$, which serves as the $L$-convex analog of the classical affine surface area. They prove, in the case when $L=RB^d$, that the quantity $as^R(K)=as^{RB^d}(K)$ is the right derivative of the volume of the $L$-convex floating body. Moreover, they show that $as^R(K)$ is rigid motion invariant and an upper semi-continuous valuation. It turns out that, in dimension $d=2$,  $\lim_{n\to\infty}\EE (f_0(K_n^{RB^2}))\, n^{-\frac 13}=C\cdot A(K)^{-\frac 13}\, as^R(K)$, where $C$ is a constant, independent of $K$. This result is completely analogous to the classical convex case and shows that a B\'ar\'any--Larman type result holds regarding the expectation of the vertex number and also for the missed area when $K$ is $C^2_+$ and $L=RB^2$ satisfy \eqref{curv-cond}. 
We conjecture that this relation holds in the general $L$-convex case and in $d$-dimensions as well. We also conjecture a series expansion for the area of $K^L(\delta)$ whose coefficients may be found with the help of Theorems~\ref{thm:main1} and \ref{thm:main3} and Corollaries~\ref{thm:main2} and \ref{thm:main4}, see concluding remarks in Section~\ref{sec:concluding-remarks}.

If the boundary of a convex body $M$ is $C^{k+1}_+$ ($k\geq 2$) smooth, then $M$ is strictly convex and $\partial M$ has a unique supporting line at each boundary point. For $u\in S^1$, let $u(M, x)$ be the outer normal unit vector at $x\in \partial M$. Under these conditions, $u(M,x)$ has a well-defined inverse $x(M,u)$ for all $u\in S^1$. 

Let $\kappa_M(x)$ denote the curvature of $\partial M$ at $x\in\partial M$. The notation $\kappa'_M(x)$ always means the derivative of $\kappa_M$ at $x\in\partial M$ with respect to arc-length such that $\partial M$ is positively oriented. In order to simplify notation, we also use $\kappa_M(u)$ and $\kappa'_M(u)$ to denote $\kappa_M(x(M,u))$ and $\kappa'_M(x(M,u))$, respectively. 

We consider two scenarios that are consistent with earlier results on this topic and in the classical convex case. First, we make the following assumption about the curvatures:
\begin{equation}\label{curv-cond}
\max_{x\in\partial L} \kappa_L (x)<1< \min_{y\in\partial K} \kappa_K(y).
\end{equation}
Note that in this case, $K$ is automatically $L$-convex; see \cite{Sch14}*{Section~3.2}.

\begin{theorem}\label{thm:main1}
    Let $k\geq 2$ be an integer, and let $K$ and $L$ be convex bodies in $\R^2$ with $C_+^{k+2}$-smooth boundaries such that \eqref{curv-cond} holds. Then,
    \[\EE \big(f_0(K_n^L)\big)=z_1(K,L)n^{\frac 13}+z_2(K,L)+\ldots+z_{k-1}(K,L)n^{-\frac{k-3}{3}}+O\big(n^{-\frac{k-2}{3}}\big)\]
    as $n\to\infty$. All coefficients $z_1(K,L),\ldots,z_k(K,L)$ can be determined explicitly. In particular
    \begin{align*}
    z_1(K,L)&=\sqrt[3]{\frac{2}{3A(K)}}\Gamma\left (\frac 53\right)\int_{S^1}\frac{(\kappa_K(u)-\kappa_L(u))^{\frac 13}}{\kappa_K(u)}\, \dx u=\sqrt[3]{\frac{2}{3A(K)}}\Gamma\left (\frac 53\right) as^L(K),\\[4pt]
    z_2(K,L)&=0,\\
    z_3(K,L)&=-\sqrt[3]{\frac{3 A(K)}{2}}\Gamma\left(\frac 73\right)\frac 15\\[4pt]
&\;\times\int_{S^1}\left(
    \frac{\kappa_L(u)(2\kappa_K^2(u)-3\kappa_K(u)\kappa_L(u)-\kappa_L^2(u))-\kappa_L''(u)}{2\kappa_K(u)\,\kappa_L(u)\,(\kappa_K(u)-\kappa_L(u))^{\frac 13}} \right.\\[4pt]
&\qquad\qquad + 
    \frac{2\kappa_L(u)(\kappa_K''(u)-\kappa_L''(u))+5\kappa_L'(u)(\kappa_K'(u)-\kappa_L'(u))}{6\kappa_K(u)\,\kappa_L(u)\,(\kappa_K(u)-\kappa_L(u))^{\frac 43}}\\[4pt]
&\left.\qquad\qquad -\frac{5(\kappa_K'(u)-\kappa_L'(u))^{2}}{9\kappa_K(u)\,(\kappa_K(u)-\kappa_L(u))^{\frac 73}}\right) \dx u.
\end{align*}
\end{theorem}

By the $L$-convex version (see \cite{FPV20}*{p. 500}) of Efron's identity \cite{efron}, Theorem~\ref{thm:main1} yields the following corollary.
\begin{corollary}\label{thm:main2}
     Under the same assumptions as in Theorem~\ref{thm:main1}, it holds that
     \[\EE \big(A(K\setminus K_n^L)\big)/A(K)=z_1(K,L)n^{-\frac 23}+\ldots+z_{k-1}(K,L)n^{-\frac{k}{3}}+O\big(n^{-\frac{k+1}{3}}\big)\]
    as $n\to\infty$.
\end{corollary}

The asymptotic formula $\lim_{n\to\infty} \EE(f_0(K_n^L))n^{-\frac 13}=z_1(K,L)$ was proved by
Fodor, Kevei and V\'\i gh \cite{FKV14} when $L=RB^2$, and by
Fodor, Papvári and V\'\i gh \cite{FPV20} in the general case. Fodor and Montenegro \cite{FM24}*{Theorem 2, p. 1410} determined $z_3(K,RB^2)$ under slightly different differentiability conditions. We note that in the formula of $z_3(K)$ on p. 1410 of \cite{FM24}, the correct coefficient of $R\kappa(x)$ in the second term of the integral should be $-3$.
For a more detailed history and results regarding expectations, variances, laws of large numbers, and central limit theorems in this probability model, we refer to \cite{BLN25+}.

The other setting we study is when $K=L$. This model significantly differs from the previous one. In order to simplify notation, we use $L_n=L_n^L$. For $u\in S^1$, let $w(u)=w_L(u)$ denote the width of $L$ in the direction $u$, that is, the distance between the two supporting lines of $L$ parallel to $u$. Let $L^*(u)\subset S^1$ be the semicircle centered at $u$, and let
\begin{equation*}
    \tilde{J}(u)=\int_{L^*(u)}\int_{L^*(u)}\frac{\left \| u_3\times u_4\right \|}{\kappa_L(u_3)\kappa_L(u_4)}\, \dx u_3\dx u_4
\end{equation*}

We prove the following expansions.

\begin{theorem}\label{thm:main3}
    Let $k\geq 2$ be an integer and $L$ a convex body in $\R^2$ with a $C^{k+2}_+$-smooth boundary. 
    Then,
    \[\EE \big(f_0(L_n)\big)=\tilde z_1(L) + \tilde z_2(L)n^{-\frac 12} + \ldots + \tilde z_{k-1}(L)n^{-\frac{k-2}{2}}+O\big(n^{-\frac{k-1}{2}}\big)\]
    as $n\to\infty$. All coefficients $\tilde z_1(L),\ldots,\tilde z_k(L)$ can be determined explicitly. In particular
    \begin{align*}
    \tilde z_1(L)&=\frac 12\int_{S^1}\frac{\tilde{J}(u)}{w^2(u)}\, \dx u,\\
    \tilde z_2(L)&=\Gamma\left (\frac 52\right )\frac{A^{\frac 12}(L)}{2}\int_{S^1}\frac{ j_2(u)}{w^{\frac 52}(u)}\, \dx u,\\
    \tilde z_3(L)&=A(L)\int_{S^1}\frac{j_3(u)}{w^3(u)}\, \dx u,
    \end{align*}
    where the quantities $j_2(u)$ and $j_3(u)$ are defined in \eqref{l2} and \eqref{l3}, respectively.
\end{theorem}

Theorem~\ref{thm:main3} yields, via the $L$-convex version of Efron's identity, the following corollary.

\begin{corollary}\label{thm:main4}
     Under the same assumptions as in Theorem~\ref{thm:main3}, it holds that
     \[\EE \big(A(L\setminus L_n)\big)/ A(L)=\tilde{z}_1(L)n^{-\frac 23}+\ldots+\tilde{z}_{k-1}(L)n^{-\frac{k}{3}}+O\big(n^{-\frac{k+1}{3}}\big)\]
    as $n\to\infty$.
\end{corollary}

The asymptotic formula $\lim_{n\to\infty} \EE (f_0(L_n))=\tilde z_1(L)$ was proved by Fodor, Kevei and V\'\i gh \cite{FKV14} in the special case when $L=RB^2$. Fodor, Papv\'ari and V\'\i gh \cite{FPV20} proved, in the general case, that $\lim_{n\to\infty} \EE (f_0(L_n))=const.$, however, the constant in \cite{FPV20}*{Formulas (1.4) and (1.5)} was incorrect, as pointed out by Marynych and Molchanov in \cite{MM22}. The correct coefficient is $\tilde{z}_1(L)$. We note that this model is an exception in the sense that Marynych and Molchanov \cite{MM22} proved an asymptotic formula for the expectation of the $f$-vector of $L_n$ in $\R^d$ for general $d$, without any smoothness assumptions. The coefficients $\tilde z_i(L)$ were determined by Fodor and Montenegro \cite{FM24} in the case when $L=RB^2$. Observe that in \cite{FM24}, if $L=RB^2$, then the only non-zero coefficients in the expansion belong to terms with integer powers of $n$ (the indexing is different in \cite{FM24}). 

We note that $\tilde z_2(L)=0$ when $L$ is centrally symmetric; see the argument in Section~\ref{sec:concluding-remarks}. Thus, in this case, the first non-zero coefficient in the expansion is $\tilde z_3(L)$. This is the reason why we carried out the calculations in order to obtain not only $\tilde z_2(L)$ but also the much more complicated $\tilde z_3(L)$.

\section{Tools}
We will use the following statement from Gruber \cite{G96}*{Lemma~1 on p. 398}, which is stated here as in Reitzner \cite{R04}, but only in the special case $d=2$.
\begin{lemma}\label{Gruber-1}
Let $M$ be a convex body in $\R^2$ with $C^{k+1}_+$ boundary for some integer $k\geq 2$. Then there exist constants $\alpha, \beta>0$ depending only on $M$ such that the following holds for every boundary point $x$ of $M$. If $x=0$ and the (unique) tangent line of $M$ at $x$ is $\R$, then there is an $\alpha$ neighborhood of $x$ in which the boundary of $M$ can be represented by a convex function $f(\sigma)$ of differentiability class $C^{k+1}$ in $\R$. Moreover, all derivatives of $f$ up to order $k+1$ are uniformly bounded by $\beta$.
\end{lemma}

Let $u\in S^1$ and assume that $M$ is in the position in Lemma~\ref{Gruber-1}. Let $f_{M,u}$ be the function that represents the boundary of $M$ in an $\alpha$-neighborhood of $x(M,u)$. Then $f_M$ is of the form
\begin{equation*}\label{f-function}
f(\sigma)=f_{M,u}(\sigma)=b_2(M,u)\sigma^2+\ldots+b_k(M, u)\sigma^k+O(\sigma^{k+1}),
\end{equation*}
where the coefficients $b_i(M)=b_i(M,u)$, $i=2,\ldots, k$ depend on $u$.

In order to determine the coefficients $b_i(M,u)$, $i=2,\ldots, k$ in terms of intrinsic quantities of $\partial M$, i.e. the curvature $\kappa_M(u)$ and its derivatives, we recall the following facts from differential geometry. Let $r(M, s)$ be the arc-length parametrization of $\partial M$ such that the orientation is counter-clockwise. Assume that $r(M,0)=x(M,u)=0$, and $r'(M, 0)$, and the unit normal vector $r''(M,0)/\kappa_M(0)=-u$ form the basis of a Cartesian coordinate system, in which we denote the coordinate along the $r'$-axis by $\sigma$, and the $r''$-axis by $\eta$, as shown in Figure~\ref{fig:sapka-sorfejtes}. For simplicity, henceforth we omit the argument of $\kappa_M$ and their derivatives. It is always understood that $\kappa_M'(s)$, $\kappa_M''(s)$, etc. mean derivatives with respect to $s$.
Then the Taylor expansion of $r(M,s)=(\sigma(s),\eta(s))$ around $s=0$ are
\begin{align}\label{eq:sigma_s}
\sigma_M(s)&=s-\kappa_M^2(0)\frac{s^3}{3!}-3\kappa_M(0)\kappa_M'(0)\frac{s^4}{4!}\notag\\
&\qquad\qquad+(\kappa_M(0)^4-3(\kappa_M(0)')^2-4\kappa_M(0)\kappa_M''(0))\frac{s^5}{5!}+\dots,\\ 
\eta_M(s)&=\kappa_M(0)\frac{s^2}{2!}+\kappa_M'(0)\frac{s^3}{3!}+
(\kappa_M''(0)-\kappa_M^3(0))\frac{s^4}{4!}\notag\\&\qquad\qquad+(-6\kappa_M^2(0)\kappa_M'(0)+\kappa_M'''(0))\frac{s^5}{5!}+\dots. \notag
\end{align}
The Taylor expansion of $r'(M,s)$ around $s=0$ is
\begin{align*}
\sigma'_M(s)&=1-\kappa_M^2(0)\frac{s^2}{2!}-3\kappa_M(0)\kappa'_M(0)\frac{s^3}{3!}\notag\\
&\qquad\qquad+(\kappa_M(0)^4-3(\kappa_M'(0))^2-4\kappa_M(0)\kappa_M''(0))\frac{s^4}{4!}+\dots,\\
\eta'_M(s)&=\kappa_M(0)s+\kappa'_M(0)\frac{s^2}{2!}+(\kappa''_M(0)-\kappa^3_M(0))\frac{s3}{3!}\notag\\
&\qquad\qquad+(\kappa_M(0)'''-6\kappa_M^2(0)\kappa_M'(0))\frac{s^4}{4!}+\dots,
\end{align*}
see, for example, \cite{D76}*{Section~1.6}.

From the equality $f_{M,u}(\sigma_M(s))=\eta_M (s)$, we can identify the coefficients $b_2(M),\ldots, b_k(M)$. In particular, if $k\geq 5$, then
\begin{align*}
b_2(M,u)&=\frac{\kappa_M(u)}{2!},\quad  b_3(M,u)=\frac{\kappa_M'(u)}{3!}, \quad 
b_4(M,u)=\frac{\kappa_M''(u)+3\kappa_M^3(u)}{4!},\\
\quad b_5(M,u)&=\frac{19\kappa_M^2(u)\kappa_M'(u)+\kappa_M'''(u)}{5!}.
\end{align*}

In order to simplify notation, we may suppress the dependence on $u$ (and thus on $x$) in the notation when we work with a fixed $u$. We will only indicate dependence when $u$ is used in the argument. We may even omit $M$ from the notation if it is clear otherwise. 

We will also use the following statement due to Gruber \cite{G96}; see also Reitzner \cite{R04}. We restate it only for $d=2$, so this is a simpler version of the original theorem:
\begin{lemma}\label{lemma:inverse}
Let 
$$\eta=\eta(\sigma)=b_m\sigma^m+\ldots+b_{k}\sigma^k+O(\sigma^{k+1})$$
for $0\leq\sigma\leq\alpha$, $2\leq m\leq k$ be a strictly increasing function. Then there are coefficients $c_1, \ldots, c_{k-m+1}$ and a constant $\gamma>0$ such that the inverse function $\sigma=\sigma(\eta)$ has the following representation
\begin{equation*}
\sigma=\sigma(\eta)=c_1\eta^{\frac 1m}+\ldots+c_{k-m+1}\eta^{\frac{k-m+1}{m}}+O\big(\eta^{\frac{k-m+2}{m}}\big)
\end{equation*}
for $0\leq \eta\leq\gamma$. The coefficients $c_1,\ldots, c_{k-m+1}$ can be determined explicitly in terms of $b_m,\ldots, b_{k}$. In particular,
\begin{itemize}
\item[i)] $c_1=\frac{1}{b_m^{\frac 1m}},$
\item[ii)] $c_2=-\frac{b_{m+1}}{mb_{m}^{\frac{m+2}{m}}}$,
\item[iii)] $c_3=-\frac{b_{m+2}}{mb_m^{\frac{m+3}{m}}}+\frac{(m+3)b_{m+1}^2}{2m^2b_m^{\frac{2m+3}{m}}}$.
\end{itemize}
\end{lemma}

The following statement is a slightly modified version of Lemma~3 on p.~402 in Gruber \cite{G96} that was adapted by Fodor and Montenegro \cite{FM24}*{Lemma~4 on p. 1417} to the current setting.

\begin{lemma}\label{lemma:beta}
Let $\beta\in\R$. There are coefficients $\gamma_1, \gamma_2,\ldots\in \R$ depending on $\beta$ which can be determined explicitly such that for a fixed $l=1, 2,\ldots$ and $0<\alpha\leq 1$,
\begin{equation*}
\int_{0}^{\alpha (n-2)} \left (1-\frac {t}{n-2}\right )^{n-2} t^\beta \, dt=\Gamma(\beta+1)+\frac{\gamma_1}{n}+\ldots+\frac{\gamma_l}{n^l}+O\left (\frac{1}{n^{l+1}}\right ), \text{ as }n\to\infty.
\end{equation*}
In particular, 
\begin{equation*}
    \gamma_1=-\frac{\Gamma(\beta+3)}{2}, \quad \gamma_2=-\frac{\Gamma(\beta+4)}{3}-2\Gamma(\beta+3).
\end{equation*}
If $\alpha$ is chosen from a closed subinterval of $(0,1]$, then the constant in $O(\cdot)$ can be chosen independently of $\alpha$.
\end{lemma}

\section{Proof of Theorem~\ref{thm:main1}}

We call $C\subset K$ an {\em $L$-cap of $K$} if $C=\cl (K\setminus (L+p))$ for some $p\in \R^2$, where $\cl(\cdot)$ is the closure of a set. It was proved in \cite{FPV20}*{Lemma 2.1, p. 502} that for an $L$-cap of $K$, $C=\cl (K\setminus (L+p))$, there exists a unique point $x\in\partial C\cap \partial K$ and $t\geq 0$ such that
$y=x-t u(K,x)\in\partial C\cap (\partial L+p)$ and $u(L+p,y)=u(K,x)$. The point $x$ is called the vertex and $t$ the height of $C$.

Let $A(u,t)$ denote the area of the $L$-cap $C(u,t)$. For each $u\in S^1$, let $t_0(u)$ be maximal such that $A(u,t_0(u))\geq 0$. It was proved in \cite{FPV20} that
\[\EE \big(f_0(K_n^L)\big)=\frac{1}{A^2(K)} \binom n2 \int_{S^1}\int_0^{t_0(u)}\left (1-\frac{A(u,t)}{A(K)}\right)^{n-2}J(u,t)\, \dx t\dx u,\]
where
\[
  J(u,t)=\left(\frac{1}{\kappa_L(u)}-\frac{1}{\kappa_K(u)}+t\right)\int_{L^*(u,t)}\int_{L^*(u,t)}\frac{\left\|u_1\times u_2\right\|}{\kappa_L(u_1)\kappa_L(u_2)}\,\dx u_1\dx u_2.
\]
Due to the at least $C^2_+$ smoothness of $\partial K$ and $\partial L$ and condition \eqref{curv-cond}, there exists a universal constant $c>0$ such that $0\leq J(u,t)\leq c$ for all $u\in S^1$ and $0<t\leq t_0(u)$.

 Let $0<\delta<A(K)$ be an arbitrary but fixed small number. Let $0<t_1$ be such that, for arbitrary $t\in [t_1,t_0(u)]$ and $u\in S^1$, $A(u,t)\geq \delta$. Then, since $t_0(u)\leq T$ for some fixed $T>0$, we obtain that
 \begin{align*}
     \int_{S^1}\int_{t_1}^{t_0(u)}\left (1-\frac{A(u,t)}{A(K)}\right)^{n-2}J(u,t)\, \dx t\dx u &\leq c\int_{S^1}\int_{t_1}^{t_0(u)}\left (1-\frac{A(u,t)}{A(K)}\right)^{n-2} \dx t\dx u\\
    &\leq 2\pi c\int_{t_1}^{T}\left (1-\frac{\delta}{A(K)}\right)^{n-2} \dx t\\
    &\leq 2\pi c T \left (1-\frac{\delta}{A(K)}\right)^{n-2},
 \end{align*}
and therefore, if $\delta$ is chosen to be sufficiently small, then
\begin{equation}\label{main-integral}
    \EE \big(f_0(K_n^L)\big)=\frac{1}{A^2(K)} \binom n2 \int_{S^1}\int_0^{t_1}\left (1-\frac{A(u,t)}{A(K)}\right)^{n-2}J(u,t)\, \dx t\dx u + O(n^{-k}).
\end{equation}

\begin{figure}[!ht]
  \centering
  \includegraphics[width=0.65\linewidth]{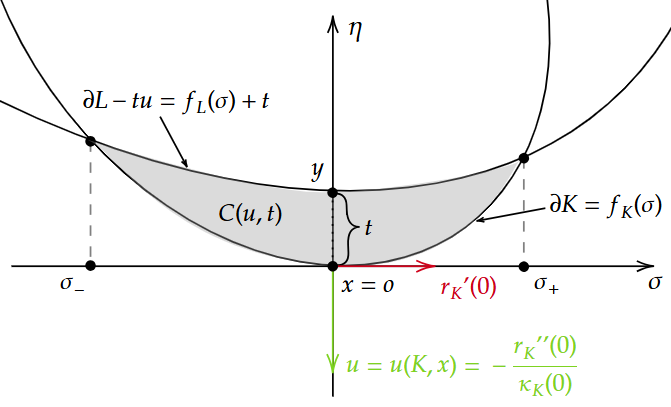}
  \caption{Position of $K$ and $L$}
  \label{fig:sapka-sorfejtes}
\end{figure}

Now, assume that $u\in S^1$, $x(K,u)=x(L,u)=0$ and $K$ and $L$ are in the position as described in Lemma~\ref{Gruber-1}, see Figure~\ref{fig:sapka-sorfejtes}. Then, for a suitable $\alpha>0$, both $\partial K$ and $\partial L$ can be represented in a common $\alpha$-neighborhood of $x$ as
\begin{align*}\label{f-function}
f_K(\sigma)&=b_2(K,u)\sigma^2+\ldots+b_k(K, u)\sigma^k+O(\sigma^{k+1}),\\
f_L(\sigma)&=b_2(L,u)\sigma^2+\ldots+b_k(L, u)\sigma^k+O(\sigma^{k+1}),
\end{align*}
where the coefficients $b_i(K)=b_i(K,u)$, and $b_i(L)=b_i(L,u)$, $i=2,\ldots, k$ depend on $u$.

Consider the $L$-cap $C(u,t)$ for sufficiently small $t>0$. Then part of the $\partial C$ belonging to $L-tu$ is represented by $t+f_L(\sigma)$
for $\sigma\in (-\alpha, \alpha)$.
Let $\sigma_-=\sigma_-(t)<0$ and  $\sigma_+=\sigma_+(t)>0$ be the $\sigma$-coordinates of the intersection points of $\partial K$ and $\partial L-tu$. Then
\[
t=t(\sigma)=(b_2(K,u)-b_2(L,u))\sigma^2+\ldots+(b_k(K,u)-b_2(L,u))\sigma^k+O(\sigma^{k+1}).
\]
For ease of notation, let $e_i:=e_i(u)=b_i(K,u)-b_i(L,u)$ for $i=1,\ldots, k$. By the condition on the curvatures $0<\kappa_L(u)<1<\kappa_K(u)$, the coefficient $e_2>0$. Thus, 
\[\sigma_+=\sigma_+(u)=c_1 t^{\frac 12}+\ldots+c_{k-1} t^{\frac{k-1}{2}}+O\big(t^{\frac k2}\big),\]
where 
\[
c_1=e_2^{-\frac 12}, \; c_2=-\frac{e_3}{2 e_2^2}, \; c_3=\frac{5 e_3^2-4 e_2 e_4}{8 e_2^{\frac 72}},\; \; c_4=-\frac{2 e_3^3-3 e_2 e_3 e_4+e_2^2 e_5}{2 e_2^5}.
\]
For the other intersection point we obtain
\[\sigma_-=\sigma_-(u)=\tilde{c}_1 t^{\frac 12}+\ldots+\tilde{c}_{k-1} t^{\frac{k-1}{2}}+O\big(t^{\frac k2}\big),\]
with $\tilde{c}_1=-c_1$, $\tilde{c}_2=c_2$, $\tilde{c}_3=-c_3$ and $\tilde{c}_4=c_4$.

Now, we are ready to calculate the area of the cap $C(u,t)$ as
\begin{align*}
    A(u,t)&=\int_{\sigma_-}^{\sigma_+} t+f_L(\sigma)-f_K(\sigma)\, \dx \sigma\notag\\
    &=\left [t\sigma-\frac{e_2}{3}\sigma^3-\ldots-\frac{e_k}{k+1}\sigma^{k+1}+O(\sigma^{k+2})\right ]_{\sigma=\sigma_-}^{\sigma_+}\notag\\
    &=a_1 t^{\frac 32}+a_2 t^2+\ldots+a_{k-1}t^{\frac{k+1}{2}}+O\big(t^{\frac{k+2}{2}}\big),
\end{align*}
where direct calculation yields that
\[
a_1=\frac 43 e_2^{-\frac 12}, \quad a_2=0, \quad a_3=\frac{5 e_3^2-4e_2e_4}{10 e_2^{\frac 72}},\quad a_4=0.
\]
Observe that, for sufficiently small $t$, $\partial A(u,t)/\partial t=\sigma_+(t)-\sigma_-(t)$. 

Let
\[
I^*(u,t)=\int_{L^*(u,t)}\int_{L^*(u,t)}\frac{\|u_1\times u_2\|}{\kappa_L(u_1)\kappa_L(u_2)}\, \dx u_1\dx u_2,
\]
and
\[
k(u,t)=\frac{1}{\kappa_L(u)}-\frac{1}{\kappa_K(u)}+t.
\]
Then, the Jacobian is $J(u,t)=k(u,t)I^*(u,t)$. Since $k(u,t)$ depends only on $u$ and $t$, we only need to examine $I^*(u,t)$. Using standard differential geometry, we get that
\begin{align*}
I^*(u,t)&=\int_{L(u,t)}\int_{L(u,t)}\|u_L(s_1)\times u_L(s_2)\|\, \dx s_1\dx s_2\\
&=\int_{s_-(t)}^{s_+(t)}\int_{s_-(t)}^{s_+(t)}\|u_L(s_1)\times u_L(s_2)\|\, \dx s_1\dx s_2,
\end{align*}
where $s_-(t)$ and $s_+(t)$ are the arc-length values on $\partial L+tu$ of the intersection points with $\partial K$.

Inverting $\sigma_L(s)$ \eqref{eq:sigma_s} by Lemma~\ref{lemma:inverse}, we obtain
\[
s=s_L(\sigma)=\sigma+\kappa_L^2\frac{\sigma^3}{3!}+ 3\kappa_L\kappa'_L\frac{\sigma^4}{4!}+
 \left(9\kappa_L^4 + 3 (\kappa _L')^2 + 4 \kappa _L \kappa _L''\right) \frac{\sigma^5}{5!}+O(\sigma^6).
\]

Then
\begin{align*}
s_+(t)&=s_L(\sigma_+(t))=g_1 t^{\frac 12}+g_2 t+g_3 t^{\frac 32}+\ldots+g_{k-1} t^{\frac {k-1}{2}}+O\big(t^{\frac{k}{2}}\big),\\
s_-(t)&=s_L(\sigma_-(t))=\tilde{g}_1 t^{\frac 12}+\tilde{g}_2 t+\tilde{g}_3 t^{\frac 32}+\ldots+\tilde{g}_{k-1} t^{\frac {k-1}{2}}+O\big(t^{\frac{k}{2}}\big),
\end{align*}
with
\begin{align*}
g_1&=e_2^{-\frac 12}, \quad g_2=-\frac{e_3}{2 e_2^2},\quad g_3=\frac{\kappa_L^2}{6 e_2^{\frac 32}}+\frac{5 e_3^2-4e_2e_4}{8 e_2^{\frac 72}},\\ g_4&=\frac{\kappa_L' \kappa _L}{8 e_2^2}-\frac{\kappa _L^2 e_3}{4 e_2^3}-\frac{2 e_3^3-3 e_2 e_3 e_4 + e_2^2 e_5}{2 e_2^5},
\end{align*}
and
\[
\tilde{g}_1=-g_1,\quad \tilde{g}_2=g_2,\quad \tilde{g}_3=-g_3,\quad \tilde{g}_4=g_4.
\]

Now, 
\begin{align*}
\|u_L(s_1)\times u_L(s_2)\|&=|\sin \angle (u_L(s_1),u_L(s_2))|=|\sin \angle (r_L'(s_1),r'_L(s_2))|\\
&=|\cos \angle (r'_L(s_1),r'^{\bot}_L(s_2))|=|\langle r'_L(s_1), r'^{\bot}_L(s_2)\rangle|\\
&=|\sigma'_L(s_1)\eta'_L(s_2)-\eta'_L(s_1)\sigma'(s_2)|,
\end{align*}
where $r'^{\bot}_L(s_2)$ is orthogonal to $r'_L(s_2)$ such that $\angle (r'_L(s_1),r'^{\bot}_L(s_2))=\frac\pi 2-\angle (r'_L(s_1),r'_L(s_2))$.

Notice that $\sigma'_L(s_1)\eta'_L(s_2)-\eta'_L(s_1)\sigma'(s_2)\geq 0$ precisely when $s_1\leq s_2$. 
Therefore,
\begin{align*}
    I^*(u,t)&=2\int_{s_1}^{s_+(t)}\int_{s_-(t)}^{s_+(t)}\sigma'_L(s_1)\eta'_L(s_2)-\eta'_L(s_1)\sigma'(s_2)\, \dx s_1\dx s_2\\
    &=l_1t^{\frac 32}+\ldots+l_{k-1}t^{\frac{k+1}{2}}+O\big(t^{\frac{k+2}{2}}\big),
\end{align*}
where, in particular,
\[
l_1=\frac{8 \kappa _L }{3 e_2^{\frac 32}},\quad 
l_2=0,\quad 
l_3=\frac{75 \kappa_L e_3^2 - 20 \kappa_L' e_2 e_3 - 60 \kappa_L e_2 e_4 + 4\left(3 \kappa_L^3+\kappa_L''\right) e_2^2}{15 e_2^{\frac 92}}.
\]
Now,
\begin{align*}
    J(u,t)&=k(u,t)I^*(u,t)=\left (\frac{1}{\kappa_L(u)}-\frac{1}{\kappa_K(u)}+t\right )I^*(u,t)\\
    &=j_1t^{\frac 32}+\ldots+j_{k-1}t^{\frac{k+1}{2}}+O\big(t^{\frac{k+2}{2}}\big),
\end{align*}
where
\begin{align*}
j_1&=\frac{8}{3 e_2^{\frac 32}}\frac{\kappa_K-\kappa_L}{\kappa_K},\quad j_2=0,\\
j_3&=\frac{8 \kappa _L}{3 e_2^{\frac 32}}+\frac{\kappa_K-\kappa_L}{\kappa_L\kappa_K}\frac{75 \kappa_L e_3^2 - 20 \kappa_L' e_2 e_3 - 60 \kappa_L e_2 e_4 + 4 \left(3 \kappa_L^3 + \kappa_L''\right) e_2^2}{15  e_2^{\frac 92}}.
\end{align*}

Next, we reparametrize the integral using cap areas. When $n$ is fixed, let $y=y(u,t)$ be the following
\[
\frac{y}{n-2}=\frac{A(u,t)}{A(K)},
\]
thus, the variable $y$ is proportional to the area of the cap $A(u,t)$. We express $t$ as a function of $\frac{y}{n-2}$ as

\begin{equation}\label{eq:t=p_i*y}
    t=p_1\left (\frac y{n-2}\right)^{\frac 23}+\ldots+ p_{k-1}\left (\frac y{n-2}\right)^{\frac k3}+O\left (\left (\frac y{n-2}\right )^{\frac{k+1}3}\right ),
\end{equation}

where,
\begin{align*}
    p_1=\left (\frac{3A(K)}{4}\right)^{\frac 23}e_2^{\frac 13},\quad p_2=0,\quad p_3=\left(\frac{3A(K)}{4}\right)^{\frac 43}\frac{4e_2e_4-5e_3^2}{20e_2^{\frac 73}},\quad p_4=0.
\end{align*}

Finally, substituting \eqref{eq:t=p_i*y} into the Jacobian, we get
\begin{equation}\label{eq:Jacobi=q_i*y}
    J\left (u,\frac y{n-2}\right )=q_1\left(\frac y{n-2}\right)+\dots+q_{k-1}\left(\frac y{n-2}\right)^{\frac{k+1}3}+O\left (\left (\frac y{n-2}\right )^{\frac{k+2}{3}}\right),
\end{equation}
with
\[
q_1=j_1p_1^{\frac 32},\quad q_2=0,\quad q_3=j_3p_1^{\frac 52}+\frac{3j_1p_3 p_1^{\frac 12}}{2}.
\]
In the expressions of the coefficients $q_1,q_3$, we use $j_1, j_3$ and $p_1,p_3$ instead of $e_i$ to simplify the notation. When evaluating the integral \eqref{main-integral}, we will substitute their values explicitly.

Now, substituting \eqref{eq:Jacobi=q_i*y} in the integral \eqref{main-integral} and using \eqref{eq:t=p_i*y}, we need to evaluate
\begin{align*}
    \EE (f_0(K_n^L))&=\frac{1}{A^2(K)}\binom n2 \int_{S^1}\int_0^{t_1}\left (1-\frac{A(u,t)}{A(K)}\right)^{n-2}J(u,t)\, \dx t\dx u+O(n^{-k})\\
    &=\frac{1}{A^2(K)}\binom n2 \frac 1{n-2}\int_{S^1}\int_0^{\tau(n-2)}\left (1-\frac{y}{n-2}\right)^{n-2}J\left (u,\frac y{n-2}\right)\\
    &\qquad\times t'\left (u,\frac y{n-2}\right)\, \dx y\dx u+O(n^{-k}),
\end{align*}
where we can obtain $t'(u,y/(n-2))$ using the formula for the derivative of inverse functions.
In the inner integral, we collect the terms according to the exponent of $\frac{y}{n-2}$, including the error term. We get
\begin{align}\notag
    &\frac{1}{2A^2(K)}\frac {n(n-1)}{n-2}\int_0^{\tau(n-2)}\left (1-\frac{y}{n-2}\right)^{n-2}J\left (u,\frac y{n-2}\right)t'\left (u,\frac y{n-2}\right)\, \dx y\\ \label{eq:integrals}
    &=\frac{1}{2A^2(K)}\frac {n(n-1)}{n-2}\int_0^{\tau(n-2)}\left (1-\frac{y}{n-2}\right)^{n-2}
    \left [v_1\left (\frac y{n-2}\right)^{\frac 23}+\ldots\right.\\ \notag
    &\left.\qquad+ v_{k-1}\left (\frac y{n-2}\right)^{\frac k3}+O\left (\left (\frac y{n-2}\right )^{\frac{k+1}3}\right )\right ]\, \dx y,
\end{align}
where
\[
v_1=\frac23 p_1q_1, \quad 
v_2=0, \quad 
v_3=\frac 43 p_3q_1+\frac 23p_1q_3.
\]
Note that the coefficients $v_i$ can also be expressed explicitly in terms of the $e_i$, and in turn, by $\kappa_K$ and $\kappa_L$ and their derivatives. In particular,
\begin{align*}
v_1&=2\cdot\sqrt[3]{\frac 23}A^{\frac 53}(K)\frac{(\kappa_K-\kappa_L)^{\frac 13}}{\kappa_K}, \qquad v_2=0,\\[6pt]
v_3&=2\cdot \sqrt[3]{\frac{3}{16}}\;A^{\frac 73}(K)\left[\frac{\kappa_L''-\kappa_L(2\kappa_K^2-3\kappa_K\kappa_L-\kappa_L^2)}{5\kappa_K\,\kappa_L\,(\kappa_K-\kappa_L)^{\frac 13}} \right.\\[4pt]
& \qquad\quad \qquad\qquad\qquad \left. - \frac{2\kappa_L(\kappa_K''-\kappa_L'')+5\kappa_L'(\kappa_K'-\kappa_L')}{15\kappa_K\,\kappa_L\,(\kappa_K-\kappa_L)^{\frac 43}}
 +\frac{2(\kappa_K'-\kappa_L')^{2}}{9\kappa_K\,(\kappa_K-\kappa_L)^{\frac 73}}\right].
\end{align*}

Then, by Lemma~\ref{lemma:beta}, we obtain that the first term in \eqref{eq:integrals} is
\begin{align*}
\frac{v_1}{2A^2(K)}&\frac {n(n-1)}{(n-2)^{\frac 53}}\int_0^{\tau(n-2)}\left (1-\frac{y}{n-2}\right)^{n-2} y^{\frac 23}\, \dx y \\
=&\sqrt[3]{\frac{2}{3A(K)}}\frac{(\kappa_K-\kappa_L)^{\frac 13}}{\kappa_K}\left(\Gamma\left (\frac 53\right )n^{\frac 13}+\left(\frac 73\Gamma\left(\frac 53\right)-\frac{\Gamma\left(\frac{10}{3}\right)}{2}\right)n^{-\frac 23}+\ldots\right).
\end{align*}

The next (non-zero) term is the following:
\begin{align*}
    \frac{v_3}{2A^2(K)}&\frac {n(n-1)}{(n-2)^{\frac 73}}\int_0^{\tau(n-2)}\left (1-\frac{y}{n-2}\right)^{n-2} y^{\frac 43}\, \dx y \\
    &=\frac{v_3}{2A^2(K)}\left(\Gamma\left (\frac 73\right )n^{-\frac 13}+\left(\frac{11\Gamma\left(\frac 73\right)}{3}-\frac{\Gamma\left(\frac{13}{3}\right)}{2}\right)n^{-\frac 43}+\ldots\right).
\end{align*}

By evaluating the $k-1$ terms in the integral \eqref{eq:integrals}, we obtain
\begin{multline*}
    \frac{1}{A^2(K)}\binom n2 \int_0^{t_1}\left (1-\frac{A(u,t)}{A(K)}\right)^{n-2}J(u,t)\, \dx t\\=w_1n^{\frac 13}+w_2+w_3n^{-\frac{1}{3}}+\ldots+w_{k-1}n^{-\frac{k-3}{3}}+O\big(n^{-\frac{k-2}{3}}\big),
\end{multline*}
where the coefficients $w_1,\ldots,w_{k-1}$ can be given explicitly. Specifically,
\begin{align*}
    w_1(u)&=\sqrt[3]{\frac{2}{3A(K)}}\Gamma\left (\frac 53\right)\frac{(\kappa_K(u)-\kappa_L(u))^{\frac 13}}{\kappa_K(u)},\\
    w_2(u)&=0,\\
    w_3(u)&=\sqrt[3]{\frac{3 A(K)}{16}}\Gamma\left(\frac 73\right)
    \left[
    \frac{\kappa_L''(u)-\kappa_L(u)(2\kappa_K^2(u)-3\kappa_K(u)\kappa_L(u)-\kappa_L^2(u))}{5\kappa_K(u)\,\kappa_L(u)\,(\kappa_K(u)-\kappa_L(u))^{\frac 13}}
    \right.\\[4pt] 
    &\qquad\qquad\qquad\qquad\quad - 
    \frac{2\kappa_L(u)(\kappa_K''(u)-\kappa_L''(u))+5\kappa_L'(u)(\kappa_K'(u)-\kappa_L'(u))}{15\kappa_K(u)\,\kappa_L(u)\,(\kappa_K(u)-\kappa_L(u))^{\frac 43}}\\[4pt] 
&\qquad\qquad\qquad\qquad\quad \left. +\frac{2(\kappa_K'(u)-\kappa_L'(u))^{2}}{9\kappa_K(u)\,(\kappa_K(u)-\kappa_L(u))^{\frac 73}}\right],
\end{align*}
where we recall that $\kappa_K$ and $\kappa_L$ are functions of $u$.

Finally, integration with respect to $u$ yields
\begin{align*}
    \EE (f_0(K_n^L))&=\int_{S^1}w_1(u)n^{\frac 13}+w_2(u)+\ldots+w_{k-1}(u)n^{-\frac{k-3}{3}}+O\big(n^{-\frac{k-2}{3}}\big)\, \dx u\\
    &=z_1(K,L)n^{\frac 13}+z_2(K,L)+\ldots+z_{k-1}(K,L)n^{-\frac{k-3}{3}}+O\big(n^{-\frac{k-2}{3}}\big),
\end{align*}
where, once again, all coefficients can be explicitly obtained by integration on $S^1$ with respect to $u$. 
This completes the proof of Theorem~\ref{thm:main1}.

\section{Proof of Theorem~\ref{thm:main3}}
Using the asymptotic formula for $\EE(f_0(L_{n}))$ proved in \cite{FPV20}*{p. 509} and a similar argument as above, we obtain that
\begin{multline}
    \EE (f_0(L_{n}))=\frac{1}{A^2(L)}\binom n2 \int_{S^1}\int_0^{t_1}\left (1-\frac{A(u,t)}{A(L)}\right)^{n-2}
    \!\!\!\!J(u,t)\, \dx t\dx u+O(n^{-k}),
    \label{eq:asymptoticK=L}
\end{multline}
where
\[J(u,t)=\int_{L(u,t)}\int_{L(u,t)}
\left\|u_3\times u_4\right\|t \,\dx u_3\dx u_4.\]

The evaluation of \eqref{eq:asymptoticK=L} differs significantly from the previous sections. 
The major difference compared to the previous case is in the nature of the cap $C(u,t)$. In the previous case, $L(u,t)$ shrank to the point $x(u)$ as $t\to 0$. However, if $K=L$, then the part of the cap $C(u,t)$ on $\partial L-ut$ "curves back", and it tends monotonically decreasingly to the arc of $\partial L$ whose endpoints are the points of tangency with the two supporting lines parallel to $u$. 
Although $L(u,t)$ is on $\partial L-tu$, we will work on its translated copy on $L$; that is, in the computation, we use $L(u,t)+tu$. 

As we work solely on the boundary of $L$, we omit $L$ from the notation if this does not cause confusion. 
For $u\in S^1$, let $x_1=x_1(u)(=x_1(L,u))$ and $x_2=x_2(u)(=x_2(L,u))$ be the points of tangency with the supporting lines of direction $u$ such that $x_1(u)$ precedes $x_2(u)$ in the positive orientation of $\partial L$. Furthermore, let $u_1=u_1(u)$ and $u_2=u_2(u)$ be the outer unit normals to $\partial L$ at $x_1(u)$ and $x_2(u)$, respectively. 
Let $r(s_1)=x_1(L,u)$ and $r(s_2)=x_2(L,u)$. For $i=1,2$, consider the Cartesian coordinate systems with origin $x_i(u)$, and axes in the directions of $r'(s_i)$ and $-u_i$. We denote the corresponding Cartesian coordinates by $\sigma_i$ and $\eta_i$ for $i=1,2$.

\begin{figure}[!ht]
    \centering
    \includegraphics[width=0.5\linewidth]{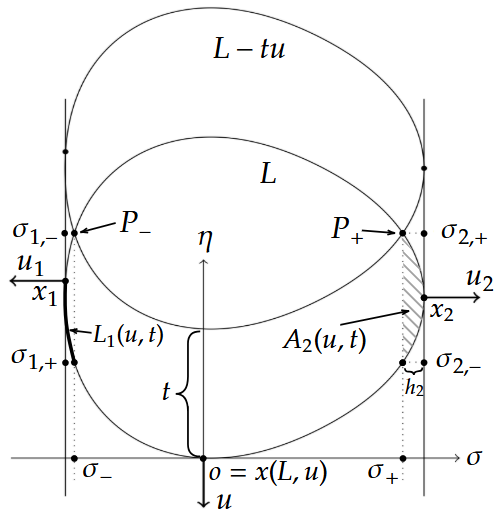}
    \caption{The $L$-cap when $K=L$}
    \label{fig:2}
\end{figure}

Then, for a sufficiently small $\alpha$, the boundary $\partial L$ can be represented by a convex function $f_i(\sigma_i)$, $i=1,2$ in the corresponding coordinate systems such that
\[
f_i(\sigma_i)=b_2(u_i)\sigma_i^2+\ldots+b_k(u_i)\sigma_i^k+O(\sigma_i^{k+1}), \quad i=1,2.
\]
From the local representation of $\partial L$ in terms of the curvature function and its derivatives, we get, similarly as before, that, for $i=1,2$, 
\begin{align*}
b_2(u_i)&=\frac{\kappa_L(u_i)}{2!}, \quad b_3(u_i)=\frac{\kappa_L'(u_i)}{3!}, \quad 
b_4(u_i)=\frac{\kappa_L''(u_i)+3\kappa_L^3(u_i)}{4!}, \\
b_5(u_i)&=\frac{19\kappa_L^2(u_i)\kappa_L'(u_i)+\kappa_L'''(u_i)}{5!}.
\end{align*}
Notice that here it is important to emphasize that we are taking the values of $\kappa_L$ and its derivatives at $u_i$, which is the unit vector obtained by rotating $u$ by $-\pi/2$ in the case of $u_1$, and $\pi/2$ in the case of $u_2$; see Figure~\ref{fig:2}.

Let $P_-=P_-(u,t)$ and $P_+=P_+(u,t)$ denote the intersection points of $\partial L$ and $\partial L-tu$ such that $P_-$ precedes $P_+$ in the positive orientation of $\partial L$. 
Let $\sigma_{-}=\sigma_{-}(u,t)=\sigma(P_-(u,t))$ and $\sigma_{+}=\sigma_{+}(u,t)=\sigma(P_+(u,t))$ be the $\sigma$-coordinates of $P_-$ and $P_+$, respectively (in the coordinate system with origin $x(u)$). 
Cut $C(u,t)$ into three pieces along the vertical lines through $P_-$ and $P_+$, and refer to the areas of the resulting pieces as $A_1(u,t)$, $A^*(u,t)$, and $A_2(u,t)$. 
First, we calculate $A_i(u,t)$ for $i=1,2$. 

Let $h_1=h_1(u,t)=\sigma_-(u,t)-\sigma(x_1(u))$ and $h_2=h_2(u,t)=\sigma(x_2(u))-\sigma_+(u,t)$, and let $\sigma_{i,-}=\sigma_{i,-}(u,t)$ and $\sigma_{i,+}=\sigma_{i,+}(u,t)$ be the $\sigma_i$ coordinates of the intersection points of the line with equation $\eta_i=h_i$ and $\partial L$. Then
\begin{align*}
\sigma_{i,+}&=\sigma_{i,+}(u,t)=c_{i,1} h_i^{\frac 12}+\ldots+c_{i,k-1} h_i^{\frac{k-1}{2}}+O\big(h_i^{\frac k2}\big),\\
\sigma_{i,-}&=\sigma_{i,-}(u,t)=\tilde{c}_{i,1} h^{\frac 12}+\ldots+\tilde{c}_{i,k-1} h^{\frac{k-1}{2}}+O\big(h_i^{\frac k2}\big),
\end{align*}
where
\begin{align*}
c_{i,1}&=b_2^{-\frac 12}(u_i),\quad 
c_{i,2}=-\frac{b_3(u_i)}{2b_2^2(u_i)},\quad 
c_{i,3}=\frac{5b_3^2(u_i)-4 b_2(u_i) b_4(u_i)}{8b_2^{\frac 72}(u_i)},\\ 
c_{i,4}&=\frac{-2b_3^3(u_i)+3b_2(u_i)b_3(u_i)b_4(u_i)-b_2^2(u_i)b_5(u_i)}{2b_2^5(u_i)},
\end{align*}
with 
$\tilde{c}_{i,1}=-c_{i,1}$, 
$\tilde{c}_{i,2}=c_{i,2}$, 
$\tilde{c}_{i,3}=-c_{i,3}$, and 
$\tilde{c}_{i,4}=c_{i,4}$.

Then $t=\sigma_{i,+}-\sigma_{i,-}$ for $i=1,2$, from which we conclude that
\[
h_i=h_i(u,t)=\mu_{i,1} t^2+\ldots+\mu_{i,k-1} t^{k}+O(t^{k+1}),
\]
where
\[
\mu_{i,1}=\frac{b_2(u_i)}{4}, \quad \mu_{i,2}=0,\quad \mu_{i,3}=-\frac{5\,b_3^2(u_i)-4\,b_2(u_i)b_4(u_i)}{64\, b_2(u_i)}, \quad \mu_{i,4}=0.
\]
Then, for the areas $A_i(u,t)$, $i=1,2$, we obtain that
\begin{align*}
A_i(u,t)&=\int_{\sigma_{i,-}}^{\sigma_{i,+}} h_i(u,t)-f_i(\sigma_i)\, \dx \sigma_i\\
&=\int_{\sigma_{i,-}}^{\sigma_{i,+}} h_i(u,t)-b_2(u_i)\sigma_i^2-\ldots-b_k(u_i)\sigma_i^k+O(\sigma_i^{k+1})\, \dx \sigma_i\\
&=a_{i,1}t^3+\ldots+a_{i,k-1}t^{k+1}+O(t^{k+2}),
\end{align*}
where
\begin{align*}
a_{i,1}=\frac{b_2(u_i)}{6}, \quad 
a_{i,2}=0,\quad
a_{i,3}=\frac{1}{80}\left (-\frac{5\,b_3^2(u_i)}{b_2(u_i)}+4\,b_4(u_i)\right ).
\end{align*}

The area $A^*(u,t)$ of the "middle part" is
\begin{align*}
    A^*(u,t)&=\big({\sigma_+(u,t)}-{\sigma_-(u,t)}\big)\,t=w(u)t-\big(h_1(u,t)+h_2(u,t)\big)\,t\\
    &=w(u)t-\big((\mu_{1,1}+\mu_{2,1})t^2+\ldots+(\mu_{1,k-1}+\mu_{2,k-1})t^{k}+O(t^{k+1})\big)\,t.
\end{align*}
Now, 
\begin{align*}
    A(u,t)=A_1(u,t)+A^*(u,t)+A_2(u,t)=a_1t+\ldots+a_{k-1}t^{k-1}+O(t^k),
\end{align*}
where
\begin{align*}
a_1&=w(u),\quad a_2=0, \quad a_3=-\frac{b_2(u_1)+b_2(u_2)}{12}, \quad a_4=0.
\end{align*}

For a fixed $n$, let $y=y(u,t)$ be the following
\[
\frac{y}{n-2}=\frac{A(u,t)}{A(L)},
\]
thus, the variable $y$ is proportional to the area of the cap $A(u,t)$. We express $t$ as a function of $\frac{y}{n-2}$ as
\begin{equation*}
    t=p_1\left (\frac y{n-2}\right)+\ldots+ p_{k-1}\left (\frac y{n-2}\right)^{k-1}+O\left (\left (\frac y{n-2}\right )^k\right ),
\end{equation*}
where, in particular,
\begin{align*}
    p_1=\frac{A(L)}{w(u)},\quad 
    p_2=0,\quad 
    p_3=\frac{b_2(u_1)+b_2(u_2)}{12\,w^4(u)}A^3(L),\quad 
    p_4=0.
\end{align*}

Furthermore, it is not difficult to see that $\frac{\partial A(u,t)}{\partial t}=\sigma_+(u,t)-\sigma_-(u,t)$ for any sufficiently small $t$. Thus,
\begin{multline*}
t'\left (\frac{y}{n-2}\right)=p_1+2p_2\left (\frac{y}{n-2}\right)+\ldots+(k-1)p_{k-1}\left (\frac{y}{n-2}\right)^{k-2}\\+O\left (\left (\frac{y}{n-2}\right)^{k-1}\right ),
\end{multline*}
where derivation is with respect to $\frac{y}{n-2}$.

Now, we turn to the calculation of the Jacobian 
\begin{align*}
J(u,t)&=\int_{L^*(u,t)}\int_{L^*(u,t)}\frac{\left\|u_3\times u_4\right\|}{\kappa_L(u_3)\kappa_L(u_4)}t\, \dx u_3\dx u_4\\
&=
t\cdot \int_{L(u,t)+tu}\int_{L(u,t)+tu}\left\|u(s_3)\times u(s_4)\right\|\, \dx s_3\dx s_4\\
&=t\cdot\int_{s_-(t)}^{s_+(t)}\int_{s_-(t)}^{s_+(t)}\|u(s_3)\times u(s_4)\|\, \dx s_3\dx s_4\\
&=t \cdot J^*(u,t),
\end{align*}
where, similarly as before, $s_-(t)$ and $s_+(t)$ are the arc-length values of $P_-(u,t)+tu$ and $P_+(u,t)+tu$ on $\partial L$, respectively.

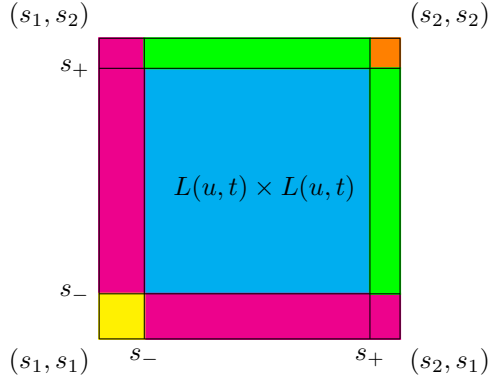
\begin{figure}[!ht]
    \centering
        \begin{tikzpicture}[scale=2]
        \filldraw [thick, color=cyan] (0.3,0.3) -- (0.3,1.8) -- (1.8,1.8) -- (1.8,0.3) -- cycle; 
        \filldraw [thick, color=magenta] (0,0.3) -- (0.3,0.3) -- (0.3,2) -- (0,2) -- cycle;
        \filldraw [thick, color=magenta] (0.3,0) -- (2,0) -- (2,0.3) -- (0.3,0.3) -- cycle;
        \filldraw [thick, color=yellow] (0,0) -- (0.3,0) -- (0.3,0.3) -- (0,0.3)  -- cycle;
        \filldraw [thick, color=orange] (1.8,1.8) -- (2,1.8) -- (2,2) -- (1.8,2)  -- cycle;
        \filldraw [thick, color=green] (1.8,0.3) -- (2,0.3) -- (2,1.8) -- (1.8,1.8)  -- cycle;
        \filldraw [thick, color=green] (0.3,1.8) -- (1.8,1.8) -- (1.8,2) -- (0.3,2)  -- cycle;
        
        \draw (0,0) -- (2,0);
        \draw (0,2) -- (2,2);
        \draw (0,0) -- (0,2);
        \draw (2,0) -- (2,2);
        \draw (0,0) node [anchor=north east] {$(s_1, s_1)$};
        \draw (2,0) node [anchor=north west] {$(s_2, s_1)$};
        \draw (0,2) node [anchor=south east] {$(s_1, s_2)$};
        \draw (2,2) node [anchor=south west] {$(s_2, s_2)$};

        \draw (0,0.3) -- (2,0.3);
        \draw (0,1.8) -- (2,1.8);
        \draw (0.3,0) -- (0.3,2);
        \draw (1.8,0) -- (1.8,2);

        \draw (0.3,0) node [anchor=north] {$s_-$};
        \draw (1.8,0) node [anchor=north] {$s_+$};
        \draw (0,0.3) node [anchor=east] {$s_-$};
        \draw (0,1.8) node [anchor=east] {$s_+$};
 
        \draw (1.1,1) node {$L(u,t)\times L(u,t)$};
        \end{tikzpicture}
    \caption{The subdivision of the integration domain of $J^*(u,t)$}
    \label{fig:domain}
\end{figure}
We subdivide the integration domain of $J^*(u,t)$ as follows.
Let $\tilde{L}(u)$ denote the arc of $\partial L$ from  $x_1(u)$ to $x_2(u)$, the two points of tangency of $\partial L$ with the supporting lines of direction $u$. This is clearly independent of $t$; it depends only on $u$. 
Furthermore, let $L_1(u,t)$ be the arc of $\partial L$ from $x_1(u)$ to $P_-(u,t)+tu$, and let $L_2(u,t)$ be the arc of $\partial L$ from $P_+(u,t)+tu$ to $x_2(u)$.
Let
\begin{align*}
\tilde{J}(u)&=\int_{\tilde{L}(u)}\int_{\tilde{L}(u)}\left \| u(s_3)\times u(s_4)\right \|\, \dx s_3\dx s_4,\\
J_i(u,t)&=\int_{L_i(u,t)}\int_{L_i(u,t)}\left \| u(s_3)\times u(s_4)\right \|\, \dx s_3\dx s_4, \quad i=1,2,\\
\overline{J}_1(u,t)&=\int_{(L(u,t)+tu)\cup L_2(u,t)}\int_{L_1(u,t)} \left \|u(s_3)\times u(s_4)\right \|\, \dx s_3\dx s_4,\\
\overline{J}_2(u,t)&=\int_{L(u,t)+tu}\int_{L_2(u,t)} \left \|u(s_3)\times u(s_4)\right \|\, \dx s_3 \dx s_4.
\end{align*}
The domains of the above integrals (taking into account symmetries) are shown in Figure~\ref{fig:domain}, where $J_1(u,t)$'s domain is yellow, $J_2(u,t)$'s is orange, $\overline{J}_1(u,t)$'s is magenta, and $\overline{J}_2(u,t)$'s is green. The domain of $J^*(u,t)$ is highlighted in cyan, and $\tilde{J}(u)$ is represented by the entire square.

With these definitions, we obtain that
\[
J^*(u,t)=\tilde{J}(u)-J_1(u,t)-J_2(u,t)-2\overline{J}_1(u,t)-2\overline{J}_2(u,t).
\]
Subsequently, we determine the above integrals one-by-one.

The quantities $J_i(u,t)$ and $i=1,2$ can be calculated similarly to how they were done in the previous Section. For ease of notation, we introduce $s_-=s(P_-+tu)$ and $s_+=s(P_++tu)$ for parametrization in the original coordinate system with the origin at $x(u)$.
Furthermore, let
\[s_{2,-}=s_+ - s_2 \quad \text{ and }\quad s_{1,+}=s_- - s_1.\]

With the notations above, we have
\begin{align*}
    J_1(u,t)&=2\int_{s_3}^{s_{-}}\int_{s_1}^{s_{-}}\sigma_1'(s_3)\eta_1'(s_4)-\eta_1'(s_3)\sigma_1'(s_4)\, \dx s_3\dx s_4\\
    &=l_{1,1}t^{3}+\ldots+l_{1,k-1}t^{k+1}+O\left (t^{k+2}\right ),
\end{align*}
and
\begin{align*}
    J_2(u,t)&=2\int_{s_3}^{s_2}\int_{s_{+}}^{s_2}\sigma_2'(s_3)\eta_2'(s_4)-\eta_2'(s_3)\sigma_2'(s_4) \,\dx s_3\dx s_4\\
    &=l_{2,1}t^{3}+\ldots+l_{2,k-1}t^{k+1}+O\left (t^{k+2}\right ),
\end{align*}
where the coefficients $l_{i,1},l_{i,2},\ldots$ can be computed explicitly. For the purposes of this paper, we present only the leading terms $l_{i,1}$, $i=1,2$. In particular,
\[ l_{i,1}=\frac{\kappa_{L}(u_i)}{24}.\]
We note that since the expansions begin at order $t^3$, the quantities $J_i(u,t)$, $i=1,2$ will not contribute to the first three terms in Theorem~\ref{thm:main3}.

Next, we calculate $\overline{J}_1(u,t)$:
This is the difficult part as it involves a "global" and a "local" part. In the "global" part, one cannot use the local expansions, the $f_i$-s. The integral $\overline{J}_1(u,t)$ can be rewritten as
\begin{align*}\overline{J}_1(u,t)&=\int_{s_-}^{s_2}\int_{0}^{s_{1,+}} \left \|u(s_3-s_1)\times u(s_4)\right \|\, \dx (s_3-s_1) \dx s_4\\
&=\int_{0}^{s_{1,+}}\int_{s_-}^{s_2}\langle u(s_3), r'(s_4)\rangle\, \dx s_4\dx s_3.
\end{align*}
For a fixed $s_3$, the outer integral is equal to the length of the projection of the arc of $\partial L$ from $P_-+tu$ to $x_2(u)$ into the line parallel to $u(s_3)$. This projection is equal to the projection of the segment $[P_-+tu,x_2(u)]$ to the same line. This can be calculated as follows. 
The vector $\overrightarrow{(P_-+tu)x_2(u)}$ is (in the original coordinate system with the origin at $x(u)$)  
\begin{align*}
&(\sigma(x_2(u)),\eta(x_2(u)))-(\sigma(x_1(u))+h_1(u,t),\eta(x_1(u))-\sigma_{1,+}(u,t))\\
&=(\sigma(x_2)-\sigma(x_1)-h_1,\eta(x_2)-\eta(x_1)+\sigma_{1,+}).
\end{align*}

The connection between the coordinate systems with origin $x_i(u)$, $i=1,2$ and the one with origin $x(u)$ is as follows. Let $Q$ be a point. Then
\begin{align*}
\sigma(Q)&=\sigma(x_1(u))+\eta_1(Q), \quad \eta(Q)=\eta(x_1(u))-\sigma_1(Q),\\
\sigma(Q)&=\sigma(x_2(u))-\eta_2(Q), \quad \eta(Q)=\eta(x_2(u))+\sigma_2(Q).
\end{align*}

Thus, since $u(s_3)=(-\eta_1'(s_3),\sigma_1'(s_3))$ in the coordinate system with origin $x_1(u)$, it is the following in the coordinate system with origin $x(u)$
\begin{align*}
u(s_3)=(\sigma(x_1)+\sigma_1'(s_3), \eta(x_1)+\eta_1'(s_3)).
\end{align*}

Thus, 
\begin{align*}
\int_{s_-}^{s_2}\langle u(s_3), r'(s_4)\rangle\, \dx s_4&=(\sigma(x_1)+\sigma_1'(s_3))(\sigma(x_2)-\sigma(x_1)-h_1)\\
&\qquad\qquad+(\eta(x_1)+\eta_1'(s_3))(\eta(x_2)-\eta(x_1)+\sigma_{1,+}).\notag
\end{align*}

Then
\begin{align*}
    \overline{J}_1(u,t)&=\int_{0}^{s_{1,+}}(\sigma(x_1)+\sigma_1'(s_3))(\sigma(x_2)-\sigma(x_1)-h_1)\\
&\qquad\qquad+(\eta(x_1)+\eta_1'(s_3))(\eta(x_2)-\eta(x_1)+\sigma_{1,+})\, \dx s_3\\
&=\overline{l}_{1,1} t^{\frac 12}+\ldots +\overline{l}_{1,k-1}t^{\frac{k-1}{2}}+O\big(t^{\frac{k}{2}}\big),
\end{align*}
where
\begin{align*}
    \overline{l}_{1,1}&=\frac{w(u)(1+\sigma(x_1))+\eta(x_1)(\eta(x_2)-\eta(x_1))}{b_2^{\frac 12}(u_1)},\\[4pt]
    \overline{l}_{1,2}&=\frac{(\eta(x_2)-\eta(x_{1}))\,\kappa_{L}(u_1)}{2\,b_{2}(u_1)}-\frac{\big(w(u) (1 + \sigma(x_1)) + \eta(x_1) (\eta(x_2) - \eta(x_1))\big)b_3(u_1)}{2b_2^2(u_1)}.
\end{align*}

Next, we calculate $\overline{J}_2(u,t)$ as follows.
\begin{align*}
    \overline{J}_2(u,t)&=\int_{s_-}^{s_+}\int_{s_{2,-}}^{0} \left \|u(s_3-s_2)\times u(s_4)\right \|\, \dx (s_3-s_2) \dx s_4\\
    &=\int_{s_{2,-}}^{0}\int_{s_-}^{s_+} \langle u(s_3), r'(s_4)\rangle\, \dx s_4\dx s_3.
\end{align*}
In this case, for a fixed $s_3$, the outer integral is the projection of the arc of $\partial L$ (or, equivalently, the segment connecting the endpoints of the arc) starting at $P_-+tu$ and ending at $P_++tu$. Then the vector $\overrightarrow{P_-P_+}$ is
\begin{align*}
\overrightarrow{P_-P_+}&=(\sigma(x_2)-h_2,\eta(x_2)+\sigma_{2,-})-(\sigma(x_1)+h_1, \eta(x_1)-\sigma_{1,+})\\
&=(\sigma(x_2)-\sigma(x_1)-h_1-h_2,
\eta(x_2)-\eta(x_1)+\sigma_{1,+}+\sigma_{2,-})\\
&=(w(u)-h_1-h_2,
\eta(x_2)-\eta(x_1)+\sigma_{1,+}+\sigma_{2,-}).
\end{align*}
Furthermore, since $u(s_3)=(-\eta_2'(s_3),\sigma_2'(s_3))$ in the coordinate system with origin $x_2(u)$, it is the following in the coordinate system with origin $x(u)$
\begin{align*}
u(s_3)=(\sigma(x_2)-\sigma_2'(s_3), \eta(x_2)-\eta_2'(s_3)).
\end{align*} 
Thus, 
\begin{align*}
\int_{s_-}^{s_+}\langle u(s_3), r'(s_4)\rangle\, \dx s_4&=
(\sigma(x_2)-\sigma_2'(s_3))(w(u)-h_1-h_2)\\
&\qquad\qquad+(\eta(x_2)-\eta_2'(s_3))(\eta(x_2)-\eta(x_1)+\sigma_{1,+}+\sigma_{2,-}).
\end{align*}
Therefore,
\begin{align*}
    \overline{J}_2(u,t)&=\int_{s_{2,-}}^{0}(\sigma(x_2)-\sigma_2'(s_3))(w(u)-h_1-h_2)\\
&\qquad\qquad+(\eta(x_2)-\eta_2'(s_3))(\eta(x_2)-\eta(x_1)+\sigma_{1,+}+\sigma_{2,-})\, \dx s_3\\
&=\overline{l}_{2,1} t^{\frac 12}+\ldots \overline{l}_{2,k-1}t^{\frac{k-1}{2}}+O\big(t^{\frac{k}{2}}\big),
\end{align*}
with
\begin{align*}
    \overline{l}_{2,1}&=\frac{w(u)(\sigma(x_2)-1)+\eta(x_2)(\eta(x_2)-\eta(x_1))}{b_2^{\frac 12}(u_2)} \\[4pt]
    \overline{l}_{2,2}&=\frac{(\eta(x_2)-\eta(x_{1}))\,\kappa_{L}(u_2)}{2\,b_{2}(u_2)}+\frac{\big(w(u) (\sigma(x_2)-1) + \eta(x_2) (\eta(x_2) - \eta(x_1))\big)b_3(u_2)}{2b_2^2(u_2)}.
\end{align*}

Thus, 
\begin{align*}
    J^*(u,t)=l_1+l_2t^{\frac 12}+l_3t+\ldots+l_{k-1}t^{\frac{k-2}{2}}+O\big(t^\frac {k-1}{2}\big),
\end{align*}
where
\begin{align}
    l_1&=\tilde{J}(u),\notag\\[4pt]
    l_2&=-2\left(
\frac{w(u)(1+\sigma(x_{1})) +\eta(x_1)(\eta(x_2)- \eta(x_1))}{b_{2}^{\frac 12}(u_1)}\right.\notag\\[4pt]
&\qquad\qquad\qquad\qquad\qquad \left. +
\frac{w(u)(\sigma(x_{2})-1) +\eta(x_2)(\eta(x_2)-\eta(x_1))}{b_{2}^{\frac 12}(u_2)}
\right),\label{l2} \\[6pt]
l_3&=\frac{\big(w(u)(1+\sigma(x_{1})) +\eta(x_1)(\eta(x_{2})-\eta(x_{1}))\big)\,b_{3}(u_1)}{b_{2}^{2}(u_1)}+\frac{(\eta(x_1)-\eta(x_{2}))\,\kappa_{L}(u_1)}{b_{2}(u_1)}
\notag\\ 
&\quad +
\frac{\big(w(u)(1-\sigma(x_{2})) + \eta(x_2)(\eta(x_1)- \eta(x_2))\big)\,b_{3}(u_2)}{b_{2}^{2}(u_2)}+\frac{(\eta(x_{1})-\eta(x_{2}))\,\kappa_{L}(u_2)}{b_{2}(u_2)}\label{l3}.
\end{align}

Therefore,
\[
J(u,t)=t\cdot J^*(u,t)=j_1t+j_2t^{\frac 32}+\ldots+j_{k-1}t^{\frac{k}{2}}+O\big(t^{\frac{k+1}{2}}\big),
\]
with $j_i=l_i$ for all $i$. In order to simplify notation, we introduce the abbreviations $j_2(u)$ and $j_3(u)$ in place of the full expressions. This not only streamlines the paper but also makes explicit the dependence on the parameter $u$.

Now, we are ready to prove Theorem~\ref{thm:main3}. We need to evaluate
\begin{align*}
    \EE (f_0(L_{n}))&=\frac{1}{A^2(L)}\binom n2 \int_{S^1}\int_0^{t_1}\left (1-\frac{A(u,t)}{A(L)}\right)^{n-2}J(u,t)\, \dx t\dx u+O(n^{-k})\\
    &=\frac{1}{A^2(L)}\binom n2 \frac 1{n-2}\int_{S^1}\int_0^{\tau(n-2)}\left (1-\frac{y}{n-2}\right)^{n-2}J\left (u,\frac y{n-2}\right)\\
    &\qquad\times t'\left (u,\frac y{n-2}\right)\, \dx y\dx u+O(n^{-k}).
\end{align*}
For the inner integral, we get
\begin{align}\notag
    &\frac{1}{2A^2(L)}\frac {n(n-1)}{n-2}\int_0^{\tau(n-2)}\left (1-\frac{y}{n-2}\right)^{n-2}J\left (u,\frac y{n-2}\right)t'\left (u,\frac y{n-2}\right)\, \dx y\\ \label{eq:integralsK=L}
    &=\frac{1}{2A^2(L)}\frac {n(n-1)}{n-2}\int_0^{\tau(n-2)}\left (1-\frac{y}{n-2}\right)^{n-2}
    \left [\tilde{v}_1\left (\frac y{n-2}\right)^{1}+\ldots\right.\\ \notag
    &\left.\qquad+ \tilde{v}_{k-1}\left (\frac y{n-2}\right)^{\frac k2}+O\left (\left (\frac y{n-2}\right )^{\frac{k+1}{2}}\right )\right ]\, \dx y,
\end{align}
where
\begin{align*}
    \tilde{v}_1&=j_1p_1^2=\tilde{J}(u)\frac{A^2(L)}{w^2(u)}
,\\[4pt]
\tilde{v}_2&=j_2 p_1^{\frac 52}=j_2(u)\left(\frac{A(L)}{w(u)}\right)^{\frac 52},
\\[4pt]
\tilde{v}_3&=j_3 p_1^3=j_3(u)\frac{A^3(L)}{w^3(u)}.
\end{align*}

Then, by Lemma~\ref{lemma:beta}, we obtain that the first term in \eqref{eq:integralsK=L} is
\begin{align*}
\frac{\tilde{v}_1}{2A^2(L)}&\frac {n(n-1)}{(n-2)^{2}}\int_0^{\tau(n-2)}\left (1-\frac{y}{n-2}\right)^{n-2} y^1\, \dx y \\
&=\frac{\tilde{J}(u)}{2w^2(u)}\frac {n(n-1)}{(n-2)^{2}}\left(\Gamma\left (2\right )-\frac{\Gamma(4)}{2} \frac 1n-\left(\Gamma(4)+\frac{\Gamma(5)}{3}+\frac{\Gamma(6)}{8}\right)\frac{1}{n^2}+\ldots\right)\\
&=\frac{\tilde{J}(u)}{2w^2(u)}\left ( 1+3n^{-1}+8n^{-2}+
\ldots\right)\left(1-3n^{-1}-29n^{-2}-\ldots\right)\\
&=\frac{\tilde{J}(u)}{2w^2(u)}\left(1-30n^{-2}+O(n^{-3})\right).
\end{align*}

The second term is 
\begin{align*}
\frac{\tilde{v}_2}{2A^2(L)}&\frac {n(n-1)}{(n-2)^{\frac 52}}
\int_0^{\tau(n-2)}\left (1-\frac{y}{n-2}\right)^{n-2} y^{\frac 32}\, \dx y \\
&=\frac{j_2 p_1^{\frac 52}}{2A^2(L)}\frac {n(n-1)}{(n-2)^{\frac 52}}
\left(\Gamma\left (\frac 52\right )-\frac 12\Gamma\left (\frac 92\right )n^{-1}
+\ldots\right)\\
&=\frac{j_2 p_1^{\frac 52}}{2A^2(L)}
\left (n^{-\frac 12}+4n^{-\frac 32}+
\ldots\right)
\left(\Gamma\left (\frac 52\right )-\frac 12\Gamma\left (\frac 92\right )n^{-1}+\ldots\right)\\
&=\frac{j_2\,p_1^{\frac 52}}{2 A^2(L)} \left(
\Gamma\left(\frac 52\right)n^{-\frac 12}
- \frac 12 \Gamma\left(\frac 92\right)n^{-\frac 32}
+ O\big(n^{-\frac 52}\big)
\right).
\end{align*}

Finally, the third term is
\begin{align*}
    \frac{\tilde{v}_3}{2A^2(L)}&\frac {n(n-1)}{(n-2)^{3}}
\int_0^{\tau(n-2)}\left (1-\frac{y}{n-2}\right)^{n-2} y^{2}\, \dx y \\
    &=\frac{j_3 p_1^3}{2A^2(L)}\frac {n(n-1)}{(n-2)^3}
\left(\Gamma(3)-\frac{\Gamma(5)}{2}\frac 1n 
-\left(\frac{\Gamma(6)}{3} + \frac{\Gamma(7)}{8}\right)\frac{1}{n^2} -\ldots\right)\\
    &=\frac{j_3 p_1^3}{2A^2(L)} 
    \left( n^{-1} + 5n^{-2}+ \ldots\right)
\left(2 - 12n^{-1} - 154n^{-2} - \ldots\right)\\
    &=\frac{j_3p_1^3}{2 A^2(L)} \left(2n^{-1} - 2n^{-2} + O(n^{-4})
\right).
\end{align*}

Thus,
\begin{align*}
    &\frac{1}{2A^2(L)}\frac {n(n-1)}{n-2}\int_0^{\tau(n-2)}\left (1-\frac{y}{n-2}\right)^{n-2}J\left (u,\frac y{n-2}\right)t'\left (u,\frac y{n-2}\right)\, \dx y\\
    &=\tilde{w}_1(u)+\tilde{w}_2(u)n^{-\frac 12}+\tilde{w}_3(u)n^{-1}+\ldots+\tilde{w}_{k-1}(u)n^{-\frac{k-2}{2}}+O\big(n^{-\frac{k-1}{2}}\big),
\end{align*}
where
\[
\tilde{w}_1(u)=\frac{\tilde{J}(u)}{2w^2(u)}, \quad \tilde{w}_2(u)=\Gamma\left (\frac 52\right )\frac{ A^{\frac 12}(L)}{2w^{\frac 52}(u)}j_2(u),\quad \tilde{w}_3(u)=\frac{A(L)}{w^3(u)}j_3(u),
\]
from which Theorem~\ref{thm:main3} follows by integration with respect to $u$.

\section{Concluding remarks}\label{sec:concluding-remarks}
We note the following phenomenon. If $L$ is centrally symmetric, then $\kappa_L(u_1)=\kappa_L(u_2)$ for all $u\in S^1$, so
\[
\tilde{w}_2(u)=-\Gamma\left (\frac 52\right )\frac{A^{\frac 12}(L)}{w^{\frac 52}(u)\kappa_L^{\frac 12}(u_1)}\big(\|x_2(u)\|-\|x_1(u)\|\big).
\]
By symmetry, 
\[
\tilde z_2(L)=\int_{S^1} \tilde{w}_2(u)\, \dx u=
-\Gamma\left (\frac 52\right )A^{\frac 12}(L)\int_{S^1} \frac{\|x_2(u)\|-\|x_1(u)\|}{w^{\frac 52}(u)\kappa_L^{\frac 12}(u_1)}\, \dx u=0,
\]
and thus, for centrally symmetric $L$, the second non-zero coefficient is $\tilde z_3(L)$. This phenomenon can also be observed when $L=R B^2$, as it was computed in \cite{FM24}.

Based on the results of Sch\"utt, Werner and Yalikun \cite{SWY25}, and on Theorems~\ref{thm:main1} and \ref{thm:main3} and Corollaries~\ref{thm:main2} and \ref{thm:main4}, the area of the $L$-convex wet part of $K$ probably has a series expansion of the form
\[
A(K^L(\delta))=f_2\delta^{\frac 23}+f_3 \delta+f_4\delta^{\frac 43}+\ldots+f_k\delta^{\frac{k-1}3}+O(\delta^{\frac k3})
\]
for $K$ and $L$ with $C^{k+2}_+$ smooth boundaries as $\delta\to 0$, cf. \cite{R01}*{Theorem~2 on p.~157}. The coefficients $f_i$ depend on $K$ and $L$.

The fact that the constants $z_i(K,L)$ for $i=1,2,3$ reduce to those obtained by Reitzner in \cite{R01} (Theorem~3 on p. 159 and the formula for $E_n(K)$ on p. 160) when $L=RB^2$ and $R\to\infty$, together with the above series expansion of $A(K^L(\delta))$, suggests that a similar relation holds for $A(K_n^L)$ as in the following equality proved by Buchta and Reitzner \cite{BR97}
\[
1-A(K_n)=\frac 43 n \int_0^{1} \left(1-\delta\right)^{n-1}A(K(\delta))\, \dx \delta,
\]
in which it is assumed that $A(K)=1$, and $K_n$ denotes a classical uniform random polygon in $K$ generated by $n$ i.i.d. random points. Using the $L$-convex version of such an equality would make it possible to determine the constants $f_i$ similarly to the approach in Reitzner \cite{R01}. We expect a very close analogy with the classical convex case here too, via the relative affine surface area.

Finally, we note that when $L$ is centrally symmetric, then $\tilde z_1(L)$ has a very transparent geometric meaning, as pointed out by Marynych and Molchanov in \cite{MM22}. Let $\Pi L$ denote the projection body of $L$, which is the plane convex body whose support function in the direction $u$ is equal to $w(u)$. Let $\Pi^\circ L$ denote the polar body of $\Pi L$. Then 
\[
\lim_{n\to\infty}\EE (f_0(L_{n}))=\tilde{z}_1(L)=\frac 12 A(\Pi L)A(\Pi^\circ L),
\]
which is proportional to the volume product of $\Pi L$.
\section*{Declarations}

\subsection*{Funding}

This research was supported by NKFIH project no.~150151, which has been implemented with the support provided by the Ministry of Culture and Innovation of Hungary from the National Research, Development and Innovation Fund, financed under the ADVANCED\_24 funding scheme.

This research was also supported by project TKP2021-NVA-09. Project no. TKP2021-NVA-09 has been implemented with the support provided by the Ministry of Innovation and Technology of Hungary from the National Research, Development and Innovation Fund, financed under the TKP2021-NVA funding scheme.

The second author was also supported by project EKÖP-24-3-SZTE-558. Project no. EKÖP-24-3-SZTE-558 has been implemented with support provided by the Ministry of Culture and Innovation of Hungary from the National Research, Development, and Innovation Fund, financed under the EKÖP-24-3 funding scheme.

\subsection*{Conflicts of interest/Competing interests}
The authors have no relevant financial or non-financial interests to disclose.

\subsection*{Data availability}
Data sharing is not applicable to this article as no data sets were generated or analyzed during the current study.

\subsection*{Author Contributions}
Both authors contributed equally in the preparation of this work.

\end{document}